\documentclass[12pt, a4paper]{article}
\title{Refined multiplicative tensor product of matrix factorizations}
\author{Yves Baudelaire Fomatati\\
\small Department of Mathematics and Statistics, University of Ottawa,\\ \small Ottawa, Ontario, Canada K1N 6N5.\\ \small yfomatat@uottawa.ca.}

\date{}
\usepackage{latexsym}
\usepackage{amsmath,amssymb}
\usepackage[square,comma,numbers,sort&compress]{natbib}
\usepackage{geometry,mathtools}
\usepackage{resizegather}

\usepackage{leqno}
\usepackage{amsfonts}

\usepackage{enumerate}

\usepackage{amsthm}
\usepackage[utf8,latin9]{inputenc}

\newcommand{\rvline}{\hspace*{-\arraycolsep}\vline\hspace*{-\arraycolsep}}

\theoremstyle{plain}
\newtheorem{remark}{Remark}[section]
\theoremstyle{plain}
\newtheorem{lemma}{Lemma}[section]
\theoremstyle{plain}
\newtheorem{proposition}{Proposition}[section]
\theoremstyle{plain}
\newtheorem{theorem}{Theorem}[section]
\theoremstyle{plain}
\newtheorem{definition}{Definition}[section]
\theoremstyle{plain}
\newtheorem{corollary}{Corollary}[section]
\theoremstyle{plain}

\theoremstyle{plain}

\newtheorem{example}{Example}[section]

\theoremstyle{plain}

\usepackage[all, 2cell]{xy}
\usepackage{txfonts}

\begin{document}
\maketitle
\begin{quote}
  \textbf{Abstract}
\end{quote}

An algorithm for matrix factorization of polynomials was proposed in \cite{fomatati2022tensor} and it was shown that this algorithm produces better results than the standard method for factoring polynomials on the class of summand-reducible polynomials.
In this paper, we improve this algorithm by refining the construction of one of its two main ingredients, namely the multiplicative tensor product $\widetilde{\otimes}$ of matrix factorizations to obtain another different bifunctorial operation denoted by $\overline{\otimes}$. We refer to $\overline{\otimes}$ as the refined multiplicative tensor product of matrix factorizations.
In fact, we observe that in the algorithm for matrix factorization of polynomials developed in \cite{fomatati2022tensor}, if we replace $\widetilde{\otimes}$ by $\overline{\otimes}$, we obtain better results on the class of summand-reducible polynomials in the sense that the refined algorithm produces matrix factors which are of smaller sizes.
\\\\
\textbf{Keywords.} Matrix factorization, tensor product, standard method for factoring polynomials, summand-reducible polynomials.\\
\textbf{Mathematics Subject Classification (2020).} 15A23, 15A69, 18A05.
\section{Introduction}

In his ground-breaking paper \cite{eisenbud1980homological} published in 1980, Eisenbud invented the concept of matrix factorization. He showed how polynomials including irreducible ones can be factorized using matrices. In light of his work, the classical factorization of a polynomial $h(x)=f(x)g(x)$ can now be viewed as a $1\times 1$ matrix factorization of $h(x)$. An example of a matrix factorization of the irreducible polynomial $t(x)=x^{2}+4$ over $\mathbb{R}[x]$ is: $$\begin{bmatrix}
    x  &  -2      \\
    2  &  x
\end{bmatrix}
\begin{bmatrix}
    x  &  2     \\
    -2  & x
\end{bmatrix}
= (x^{2} + 4)\begin{bmatrix}
    1  &  0      \\
    0  &  1
\end{bmatrix}
=tI_{2} $$
We say that
$
(\begin{bmatrix}
    x  &  -2     \\
    2  &  x
\end{bmatrix},
\begin{bmatrix}
    x  &  2      \\
    -2  & x
\end{bmatrix})
$
 is a $2\times 2$ matrix factorization of $t$.\\

  There exists a standard technique for factoring polynomials using matrices (cf. subsection \ref{sec: the std method}). One noticeable downside of this algorithm is that for each monomial that is added to the polynomial, the size of the matrix factors doubles. As will be seen below (subsection \ref{sec: the std method}), polynomials with $n$ monomials have matrix factors of size $2^{n-1}$ (when $n=10$, we have matrix factors of size $2^{10-1}= 512$).  \cite{crisler2016matrix} gives an elementary way to produce matrix factors that can be found in \cite{buchweitz1987cohen} on the class of sums of square polynomials $f=x_{1}^{2}+\cdots + x_{n}^{2}$ over the reals, with $n\leq 8$. The matrix factors they obtain are smaller in size than the ones one obtains with the standard technique. This technique has recently been improved in \cite{fomatati2022tensor} on the class of \textit{summand-reducible polynomials} (cf. definition \ref{defn summand reducible polynomials}).
  Our main objective in this paper is to refine the algorithm proposed in \cite{fomatati2022tensor} for the foregoing class of polynomials. More precisely, we want to produce matrix factors of a given polynomial which are as small as possible while maintaining the property that the entries are monomials. To this end, we will first assume (cf. subsection \ref{sec: the std method}) that before applying the standard method to find matrix factors of a polynomial, it has to be written in its expanded form. This is important because when a polynomial is not written in its expanded form, the matrix factors one obtains contain entries which are more complex (sum and products of monomials) than the ones (monomials) obtained when the method is applied on a polynomial in its expanded form (cf. examples \ref{exple: good matrix facto of g} and \ref{exple: bad matrix facto of g} ). Since a polynomial is made up of sums and products of monomials, the procedure developed in \cite{fomatati2022tensor} uses two main ingredients which are functorial operations: The Yoshino tensor product of matrix factorizations (cf. subsection \ref{subsec: reduced mult tens prodt}) with the ability to produce a matrix factorization of the sum of two polynomials from their respective matrix factorizations
   and the multiplicative tensor product of matrix factorizations (cf. definition \ref{defn of the multiplicative tensor product}) with the ability to produce a matrix factorization of the product of two polynomials from their respective matrix factorizations. In this paper, we refine the construction of the multiplicative tensor product of matrix factorizations $\widetilde{\otimes}$ and obtain a new bifunctorial operation denoted $\overline{\otimes}$. We refer to  $\overline{\otimes}$ as the \textit{refined multiplicative tensor product of matrix factorizations}. In fact, $\widetilde{\otimes}$ is simply a direct sum of two copies of $\overline{\otimes}$.
    We then replace $\widetilde{\otimes}$ by its new version $\overline{\otimes}$ in the algorithm proposed in \cite{fomatati2022tensor} to obtain a refined algorithm which yields smaller matrix factors than the ones yielded by the procedure in \cite{fomatati2022tensor} on the class of \textit{summand-reducible polynomials}.

 In the sequel, except otherwise stated, $K$ is a field and $K[[x]]$ (resp. $K[[y]]$) is the formal power series ring in the variables $x=x_{1},\cdots,x_{r}$ (resp. $y=y_{1},\cdots,y_{s}$). Let $f\in K[[x]]$ and $g\in K[[y]]$ be nonzero noninvertible
 elements.\\
 Eisenbud also found out that matrix factorizations of a power series $f\in K[[x]]$ are closely related to the homological properties of modules over quotient rings $K[[x]]/(f)$. He proved that all maximal Cohen-Macaulay modules (MCM modules) without free summands are described by matrix factorizations.
 See \cite{leuschke2012cohen} and \cite{huneke2002two} for more background on MCM modules.
\\  Yoshino \cite{yoshino1998tensor} found a way to relate MCM modules over $K[[x]]/(f)$ and over $K[[y]]/(g)$ with MCM modules over $K[[x,y]]/(f+g)$. In fact, he constructed a tensor product denoted $\widehat{\otimes}$ which is such that if $X$ is a matrix factorization of $f\in K[[x]]$ and $Y$ is a matrix factorization of $g\in K[[y]]$, then $X\widehat{\otimes}Y$ is a matrix factorization of $f+g\in K[[x,y]]$. \\
In \cite{fomatati2022tensor},
 without resorting to homological methods, the author constructed a bifunctorial operation denoted by $\widetilde{\otimes}$ which is such that
$X \widetilde{\otimes} Y$ is a matrix factorizations of the product $fg\in K[[x,y]]$. \\
In this paper, we refine the construction of $\widetilde{\otimes}$ and obtain another bifunctorial operation which unlike $\widetilde{\otimes}$ has no variants as will be discussed in Remark \ref{remark: comparison between the two multiplicative tensor prdts} .

 Thus, our first main result is the following:
 \begin{quote}
   \textbf{Theorem A.}\\
1.
Let $f\in K[[x_{1},...,x_{r}]]$ and $g\in K[[y_{1},...,y_{s}]]$ be nonzero elements.
If $X$ (resp. $Y$) is a matrix factorization of $f$ (resp. $g$). Then, there is a tensor product $\overline{\otimes}$ of matrix factorizations which produces a matrix factorization $X \overline{\otimes} Y$ of the product $fg\in K[[x_{1},...,x_{r},y_{1},...,y_{s}]]$ such that the size of each matrix factor of $X \overline{\otimes} Y$ is one half the size of matrix factors of  $ X\widetilde{\otimes} Y$. $\overline{\otimes}$ is called the \textit{refined multiplicative tensor product of matrix factorizations}. \\
2. The refined multiplicative tensor product $(-) \overline{\otimes} (-):MF(K[[x]],f)\times MF(K[[y]],g)\rightarrow MF(K[[x,y]],fg)$ is a bifunctor.
 \end{quote}

 We use the newly defined operation $\overline{\otimes}$  together with the existing Yoshino tensor product $\widehat{\otimes}$, to improve the algorithm for matrix factorization of polynomials proposed in \cite{fomatati2022tensor} on the class of \textit{summand-reducible polynomials} (cf. Definition \ref{defn summand reducible polynomials}).\\

Our second main result is stated as follows:
 \begin{quote}
   \textbf{Theorem B.}
Let $f=t_{1}+\cdots + t_{s}+ g_{11}\cdots g_{1m_{1}} + \cdots + g_{l1}\cdots g_{lm_{l}}$ be a \textit{summand-reducible polynomial}. Let $p_{ji}$ be the number of monomials in $g_{ji}$. Then
there is an improved version of the standard method for factoring $f$ which produces factorizations of size $$2^{\prod_{i=1}^{m_{1}}p_{1i} + \cdots + \prod_{i=1}^{m_{l}}p_{li} - (\sum_{i=1}^{m_{1}}p_{1i} + \cdots + \sum_{i=1}^{m_{l}}p_{li} -  \sum_{j=1}^{l}m_{j} + l)}$$ times smaller than the size one would normally obtain with the standard method.\\
 \end{quote}
 \textbf{Nota Bene:} As we will explain later (see Corollary \ref{Cor: comparison between the two algorithms} ), this theorem is equivalent to saying that for a given summand-reducible polynomial $f$, the algorithm we construct in this paper produces matrix factors of $f$ whose size is $2^{\sum_{j=1}^{l}m_{j}-l}$  times smaller than the size one obtains with the improved algorithm presented in section 4 of \cite{fomatati2022tensor}.\\

Due to its applications in pure mathematics and physics, the study of matrix factorizations evolved rapidly over the last four decades (\cite{eisenbud1980homological},\cite{kapustin2003topological}, \cite{kapustin2004d}, \cite{orlov2009derived}, \cite{orlov2012matrix}, \cite{orlov2006triangulated}, \cite{dao2013vanishing}, \cite{herzog1991linear}, \cite{buchweitz1987cohen}, \cite{khovanov2008matrix},
\cite{aspinwall2012quivers}, \cite{avramov1989modules}, \cite{yu2013geometric}, \cite{carqueville2016adjunctions}, \cite{crisler2016matrix}, \cite{fomatati2019multiplicative}) as discussed in the introduction of \cite{fomatati2022tensor}.\\
The rest of this paper is organized as follows: In section $2$, we recall the definition of the Yoshino tensor product of matrix factorizations. Next, we recall the definition of the multiplicative tensor product of matrix factorizations.
Furthermore, we define the refined version of the multiplicative tensor product of matrix factorizations.
Our theorem A is also stated and proved here.
In section $3$, properties of the refined version of the multiplicative tensor product of matrix factorizations  are discussed. In section $4$, the definition of the class of summand-reducible polynomials is recalled. Our theorem B is also stated and proved here. Examples are provided to illustrate this result.

\section{Refined multiplicative tensor product of matrix factorizations and its functoriality}
In this section, we first recall the definitions of Yoshino's tensor product of matrix factorizations denoted $\widehat{\otimes}$. Next, we recall the definition of the multiplicative tensor product of matrix factorizations $\widetilde{\otimes}$ and then refine it to obtain a new version denoted by $\overline{\otimes}$. We prove that it is a bifunctorial operation and give some examples.\\
Under this section, unless otherwise stated, $R=K[[x]]$ and $S=K[[y]]$ where $x=x_{1},...,x_{r}$ and $y=y_{1},...,y_{s}$.

\subsection{A refined version of the multiplicative tensor product of matrix factorizations } \label{subsec: reduced mult tens prodt}
In this subsection, we refine the definition of the multiplicative tensor product of matrix factorizations to obtain a bifunctorial operation denoted $\overline{\otimes}$ which is referred to as the \textit{refined multiplicative tensor product of matrix factorizations}.  $\widetilde{\otimes}$ is simply a direct sum of two copies of $\overline{\otimes}$. The functoriality of this new operation will be proved in subsection \ref{sec: functoriality of reduced tens prod}.
\\
We first recall the definitions of matrix factorization of a power series, Yoshino tensor product and multiplicative tensor product.
\begin{definition} \cite{yoshino1998tensor}  \label{defn matrix facto of polyn}   \\
An $n\times n$ \textbf{matrix factorization} of a power series $f\in \;R$ is a pair of $n$ $\times$ $n$ matrices $(\phi,\psi)$ such that
$\phi\psi=\psi\phi=fI_{n}$, where $I_{n}$ is the $n \times n$ identity matrix and the coefficients of $\phi$ and of $\psi$ are taken from $R$.
\end{definition}
The category of matrix factorizations of a power series $f\in R=K[[x]]:=K[[x_{1},\cdots,x_{n}]]$ denoted by $MF(R,f)$ or $MF_{R}(f)$ was defined in $\S 1$ of \cite{yoshino1998tensor}.
Details on this category are found in chapter 2 of \cite{fomatati2019multiplicative}.
\\
\begin{definition} \cite{yoshino1998tensor} \label{defn Yoshino tensor prodt} \textbf{Yoshino tensor product of matrix factorizations}\\
Let $X=(\phi,\psi)$ be an $n\times n$ matrix factorization of $f\in R$  and $X'=(\phi',\psi')$ an $m\times m$ matrix factorization of $g\in S$. These matrices can be considered as matrices over $L=K[[x,y]]$ and the \textbf{tensor product} $X\widehat{\otimes} X'$ is given by\\
(\(
\begin{bmatrix}
    \phi\otimes 1_{m}  &  1_{n}\otimes \phi'      \\
   -1_{n}\otimes \psi'  &  \psi\otimes 1_{m}
\end{bmatrix}
,
\begin{bmatrix}
    \psi\otimes 1_{m}  &  -1_{n}\otimes \phi'      \\
    1_{n}\otimes \psi'  &  \phi\otimes 1_{m}
\end{bmatrix}
\))\\
where each component is an endomorphism on $L^{n}\otimes L^{m}$.
\end{definition}
$X\widehat{\otimes} X'$ is an object of $MF_{L}(f+g)$ of size $2nm$ as proved in Lemma 2.1 of \cite{fomatati2019multiplicative}.\\
Functorial properties of $\widehat{\otimes}$ were proved in \cite{yoshino1998tensor}. Some of them include commutativity, distributivity and associativity.
\begin{definition} \cite{fomatati2022tensor} \label{defn of the multiplicative tensor product}
Let $X=(\phi,\psi)$ be a matrix factorization of $f\in K[[x]]$ of size $n$ and let $X'=(\phi',\psi')$ be a matrix factorization of $g\in K[[y]]$ of size $m$. Thus, $\phi,\psi,\phi' \,and \,\psi'$ can be considered as matrices over $L=K[[x,y]]$ and the \textbf{multiplicative tensor product} $X\widetilde{\otimes} X'$ is given by \\\\
\[((\phi\otimes\phi')\oplus (\phi\otimes\phi'), (\psi\otimes\psi')\oplus (\psi\otimes\psi'))=(
\begin{bmatrix}
    \phi\otimes\phi'  &          0      \\
    0                  &\phi\otimes\phi'
\end{bmatrix},
\begin{bmatrix}
    \psi\otimes\psi'  &    0      \\
    0                 &  \psi\otimes\psi'
\end{bmatrix}
)\]
\\\\
where each component is an endomorphism on $L^{n}\otimes_{L} L^{m}$.
\end{definition}
Functorial properties of $\widetilde{\otimes}$ were proved in \cite{fomatati2022tensor}. Some of them include commutativity, distributivity and associativity.
\begin{definition}   \label{defn of the reduced mult tens prodt}
Let $X=(\phi,\psi)$ be a matrix factorization of $f\in K[[x]]$ of size $n$ and let $X'=(\phi',\psi')$ be a matrix factorization of $g\in K[[y]]$ of size $m$. Thus, $\phi,\psi,\phi' \,and \,\psi'$ can be considered as matrices over $L=K[[x,y]]$. The \textbf{refined multiplicative tensor product} $X\overline{\otimes} X'$ is given by \\\\
$X\overline{\otimes} X'$ = $(\phi,\psi)\overline{\otimes} (\phi',\psi')=([\phi\otimes\phi'], [\psi\otimes\psi'])$
\\\\
where each component is an endomorphism on $L^{n}\otimes_{L} L^{m}$.
\end{definition}

\begin{remark} \textbf{Comparison between $\widetilde{\otimes}$ and $\overline{\otimes}$} \label{remark: comparison between the two multiplicative tensor prdts}
\begin{itemize}
  \item Observe that unlike with $\widetilde{\otimes}$ (cf. definition \ref{defn of the multiplicative tensor product}), the entries of the matrices we have in the foregoing definition cannot be rotated to obtain a variant for $\overline{\otimes}$.
  \item Though $\widetilde{\otimes}$ and $\overline{\otimes}$ can both be used to find a matrix factorization of the product of two polynomials from their respective matrix factorizations, there is a conspicuous difference between them in the way they are defined. But most importantly, there is a striking difference in their applications. $\widetilde{\otimes}$ was used in \cite{fomatati2021some} (cf. section 4.2) to give an example of a semi-unital semi-monoidal category. Observe that $\widetilde{\otimes}$ was used to construct objects of the foregoing category but this is not possible with $\overline{\otimes}$.
      On the other hand, as we shall see in section \ref{sec: the refined algorithm}, $\overline{\otimes}$ helps to obtain smaller matrix factors of polynomials as compared to $\widetilde{\otimes}$.
\end{itemize}

\end{remark}
\begin{remark}
The proof of Theorem 4.2 of \cite{fomatati2022tensor} gives an algorithm that yields small matrix factors (as compared to the standard method) for a summand-reducible polynomial (cf. Definition \ref{defn summand reducible polynomials}). The main ingredients of that algorithm are the bifunctorial operations $\widehat{\otimes}$ and $\widetilde{\otimes}$. The main result of this paper stipulates that if in the proof of Theorem 4.2 of \cite{fomatati2022tensor}, $\widetilde{\otimes}$ is replaced by $\overline{\otimes}$, then we will obtain better results (cf. Theorem \ref{thm refined algo for summand-red polyn} and Corollary \ref{Cor: comparison between the two algorithms}) in the sense that matrix factors will even be smaller.\\
\textbf{N.B.} Now, though $\widetilde{\otimes}$ is simply a direct sum of two copies of $\overline{\otimes}$, the refined algorithm we present in this paper produces matrix factors which are not always just half the size of matrix factors obtained from the previous algorithm, instead the size is reduced by a factor that depends on the number of terms in the summand-reducible polynomial as proved in Corollary \ref{Cor: comparison between the two algorithms}.
\end{remark}
\begin{lemma} \label{lemma size of X overline tensor Y}
Let $X=(\phi,\psi)$ be an $n\times n$ matrix factorization of $f\in K[[x]]$ and let $X'=(\phi',\psi')$ be an $m\times m$ matrix factorization of $g\in K[[y]]$. Then,
$X\overline{\otimes} X'$ is an object of $MF(K[[x,y]], fg)$ of size $nm$.
\end{lemma}
\begin{proof}

We have:\\\\
$[\phi\otimes\phi'] [\psi\otimes\psi']$
\\
$=[\phi\psi\otimes\phi'\psi'],\,\,\,by\,the\,mixed\,product\,property.$
\\
$=[f1_{n}\otimes g1_{m}],\,\,\,since\,\phi\psi=f1_{n}\,\,and\,\phi'\psi'=g1_{m}$\\
$=fg[1_{n}\otimes 1_{m}].$\\
$ =fg\cdot 1_{nm}$

So, $X\overline{\otimes} X'$ is an object of $MF(fg)$ of size $nm$ as claimed.
\end{proof}
In the following example, we exhibit a matrix factorization of the polynomial\\ $g=xy+xz^{2}+yz^{2}$ without showing how it is obtained. Details showing how such matrix factorizations are obtained are discussed in section \ref{sec: the std method} where we recall the standard method for matrix factorization of polynomials. In that section, we will also observe that this method has variants.
\begin{example} 
A straightforward computation shows that a $2 \times 2$ matrix factorization of $f=x^{2}+4$ is:\\
$
X=(\begin{bmatrix}
     x  &  -2     \\
    2  &  x
\end{bmatrix},
\begin{bmatrix}
    x  &  2      \\
    -2  & x
\end{bmatrix})=(\phi_{X},\psi_{X}).
$\\\\
Since $
(\begin{bmatrix}
     z^{2}  &  y     \\
     x  &  -x-y
\end{bmatrix}
\begin{bmatrix}
    x+y  &  y      \\
    x  &  -z^{2}
\end{bmatrix})=(xy+xz^{2}+yz^{2})\begin{bmatrix}
    1  &  0      \\
    0  &  1
\end{bmatrix},
$\\
 a $2\times 2$ matrix factorization of $g=xy+xz^{2}+yz^{2}$ is:\\
 $Y=
(\begin{bmatrix}
     z^{2}  &  y     \\
     x  &  -x-y
\end{bmatrix},
\begin{bmatrix}
    x+y  &  y      \\
    x  &  -z^{2}
\end{bmatrix})=(\phi_{Y},\psi_{Y}).
$\\\\
\[ X\overline{\otimes} Y=(\begin{bmatrix}
    \phi_{X}\otimes \phi_{Y}
\end{bmatrix},
\begin{bmatrix}
    \psi_{X}\otimes \psi_{Y}
\end{bmatrix})
\]\\
\begin{gather*}
  \setlength{\arraycolsep}{1.0\arraycolsep}
  \text{\footnotesize$\displaystyle
  i.e.,\, X\overline{\otimes} Y =(
    \begin{bmatrix}
    xz^{2} & xy & -2z^{2}  &  -2y      \\
    x^{2}  &  -x^{2}-xy & -2x & 2(x+y) \\
    2z^{2}  & 2y  &   xz^{2}& xy \\
    2x  &-2(x+y)  &   x^{2}   &  -x^{2}-xy
\end{bmatrix},
\begin{bmatrix}
    x^{2}+xy & xy & 2(x+y)  &   2y      \\
    x^{2}  &  -xz^{2} & 2x &  -2z^{2} \\
    -2(x+y)  & -2y    &  x^{2}+xy & xy \\
    2x  & 2z^{2}    &   x^{2}   &  -xz^{2}
\end{bmatrix} )$}
\end{gather*}
\end{example}
\subsection{Functoriality of $\overline{\otimes}$ } \label{sec: functoriality of reduced tens prod}
This subsection is entirely devoted to the discussion of the bifunctoriality of $\overline{\otimes}$. \\

  \textbf{Setting the stage:}
Let $X_{f}=(\phi,\psi)$, $X'_{f}=(\phi',\psi')$ and $X_{f}''=(\phi'',\psi'')$ be objects of $MF(K[[x]],f)$ respectively of sizes $n, n'$ and $n''$. Let $X_{g}=(\sigma,\rho)$, $X_{g}'=(\sigma',\rho')$ and $X_{g}''=(\sigma'',\rho'')$ be objects of $MF(K[[y]],g)$ respectively of sizes $m, m'$ and $m''$.

\begin{definition}\label{defn zeta is a bifuntor}
For morphisms $\zeta_{f}=(\alpha_{f}, \beta_{f}): X_{f}=(\phi,\psi) \rightarrow X_{f}'=(\phi',\psi')$  and $\zeta_{g}=(\alpha_{g}, \beta_{g}): X_{g}=(\sigma,\rho) \rightarrow X_{g}'=(\sigma',\rho')$ respectively in $MF(K[[x]],f)$ and $MF(K[[y]],g)$,\\ we define $\zeta_{f}\overline{\otimes} \zeta_{g}: X_{f}\overline{\otimes}X_{g} =(\phi,\psi)\overline{\otimes} (\sigma,\rho)\rightarrow X_{f}'\overline{\otimes}X_{g}' =(\phi',\psi')\overline{\otimes} (\sigma',\rho')$
by
$$([\alpha_{f}\otimes \alpha_{g}],[\beta_{f}\otimes \beta_{g}])
 $$
\end{definition}

\begin{lemma}\label{zeta is a morphism in both arguments}

  $\zeta_{f}\overline{\otimes} \zeta_{g}$: $X_{f}\overline{\otimes}X_{g} =(\phi,\psi)\overline{\otimes} (\sigma,\rho)\rightarrow X_{f}'\overline{\otimes}X_{g}' =(\phi',\psi')\overline{\otimes} (\sigma',\rho')$ is a morphism in $MF(K[[x,y]],fg)$.

\end{lemma}
\begin{proof}
 We need to show that the following diagram commutes:

 $$\xymatrix@ R=0.8in @ C=1.10in{K[[x,y]]^{nm} \ar[r]^{[\psi\otimes\rho ]} \ar[d]_{[ \alpha_{f} \otimes \alpha_{g}]} &
K[[x,y]]^{nm} \ar[d]^{[\beta_{f} \otimes \beta_{g}]} \ar[r]^{[\phi\otimes\sigma]} & K[[x,y]]^{nm}\ar[d]^{[\alpha_{f} \otimes \alpha_{g}]}\\
K[[x,y]]^{n'm'} \ar[r]_{[\psi'\otimes\rho']} & K[[x,y]]^{n'm'}\ar[r]_{[\phi'\otimes\sigma']} & K[[x,y]]^{n'm'}}$$

 viz. both the left and the right squares in the foregoing diagram commute.\\
 $\bullet$ The commutativity of the right square and the left square are respectively expressed by the following equalities:\\
 $ [\alpha_{f} \otimes \alpha_{g} ][\phi\otimes\sigma]=[\phi'\otimes\sigma'][\beta_{f} \otimes \beta_{g}]$
and\\
 $[\beta_{f} \otimes \beta_{g}][\psi\otimes\rho]=[\psi'\otimes\rho'][\alpha_{f}\otimes \alpha_{g}]$.

i.e., all we need to show is the pair of equalities:
\\
$$\begin{cases}
\alpha_{f}\phi\otimes \alpha_{g}\sigma=\phi'\beta_{f}\otimes\sigma'\beta_{g} \cdots (1)\\
\beta_{f}\psi \otimes \beta_{g}\rho= \psi'\alpha_{f}\otimes\rho'\alpha_{g} \cdots (2)
\end{cases}$$

Now by hypothesis, $\zeta_{f}=(\alpha_{f},\beta_{f}): X_{f}=(\phi,\psi) \rightarrow X_{f}'=(\phi',\psi')$ and $\zeta_{g}=(\alpha_{g},\beta_{g}): X_{g}=(\sigma,\rho) \rightarrow X_{g}'=(\sigma',\rho')$ are morphisms, meaning that the following diagrams commute
$$\xymatrix@ R=0.6in @ C=.75in{K[[x]]^{n} \ar[r]^{\psi} \ar[d]_{\alpha_{f}} &
K[[x]]^{n} \ar[d]^{\beta_{f}} \ar[r]^{\phi} & K[[x]]^{n}\ar[d]^{\alpha_{f}}\\
K[[x]]^{n'} \ar[r]^{\psi'} & K[[x]]^{n'}\ar[r]^{\phi'} & K[[x]]^{n'}}$$

and \\
$$\xymatrix@ R=0.6in @ C=.75in{K[[y]]^{m} \ar[r]^{\rho} \ar[d]_{\alpha_{g}} &
K[[y]]^{m} \ar[d]^{\beta_{g}} \ar[r]^{\sigma} & K[[y]]^{m}\ar[d]^{\alpha_{g}}\\
K[[y]]^{m'} \ar[r]^{\rho'} & K[[y]]^{m'}\ar[r]^{\sigma'} & K[[y]]^{m'}}$$

  That is,
$$\begin{cases}
 \alpha_{f}\phi=\phi'\beta_{f} \cdots (i) \\
 \psi'\alpha_{f}= \beta_{f}\psi \cdots (ii)
\end{cases}$$
and \\
$$\begin{cases}
 \alpha_{g}\sigma=\sigma'\beta_{g} \cdots (i') \\
 \rho'\alpha_{g}= \beta_{g}\rho \cdots (ii')
\end{cases}$$
Now considering $(i)$ and $(i')$, we immediately see that equality $(1)$ holds. Similarly, $(ii)$ and $(ii')$ yield $(2)$. \\

So, $\zeta_{f} \overline{\otimes} \zeta_{g}$ is a morphism in $MF(K[[x,y]], fg)$.

\end{proof}
We can now state the following result.
\begin{theorem} \label{zeta is a bifunctor}
\begin{enumerate}
\item Let $X$ be a matrix factorization of $f\in K[[x]]$ of size $n$ and let $Y$ be a matrix factorization of $g\in K[[y]]$ of size $m$. Then, there is a tensor product $\overline{\otimes}$ of matrix factorizations which produces a matrix factorization $X \overline{\otimes} Y$ of the product $fg\in K[[x_{1},...,x_{r},y_{1},...,y_{s}]]$ which is of size $nm$. $\overline{\otimes}$ is called the refined multiplicative tensor product of matrix factorizations.
\item The refined multiplicative tensor product $(-) \overline{\otimes} (-):MF(K[[x]],f)\times MF(K[[y]],g)\rightarrow MF(K[[x,y]],fg)$ is a bifunctor.
\end{enumerate}
\end{theorem}
\begin{proof}
\begin{enumerate}
  \item This is exactly what we proved above in subsection \ref{subsec: reduced mult tens prodt}.
  \item We show that $\overline{\otimes}$ is a bifunctor.\\
In order to ease our computations, let's write $F=(-)\overline{\otimes} (-)$. We show that $F$ is a bifunctor.\\
 We have:

$\;\; (-)\overline{\otimes} (-):\;\;\;\;\;\;\;\;\;\;\;\;\; MF(f)\times MF(g)\,\,\,\,\,\,\,\,\,\,\longrightarrow \,\,\,\,\,\,\,\,\,\,\,\,\,\,\,\,MF(fg)$
$$\xymatrix @ R=0.4in @ C=.33in
{&(X_{f}\;\;\;\;\;\;{,} \ar[d]_{\zeta_{f}}\ar @{}[dr]& X_{g}) \ar[rrr]^-{} \ar[d]_{\zeta_{g}}
\ar @{}[dr] &&& X_{f}\overline{\otimes} X_{g} \ar[d]^{\zeta_{f}\overline{\otimes} \zeta_{g}:= (\alpha,\beta)}\\
&(X_{f}'\;\;\;\;\;\; {,}& X_{g}') &\ar[r]&& X_{f}'\overline{\otimes} X_{g}'}$$
$$\xymatrix @ R=0.4in @ C=.3in
{&\ar[d]_{\zeta_{f}'}\ar @{}[dr]& \ar[d]_{\zeta_{g}'}
\ar @{}[dr] &&&\ar[d]^{\zeta_{f}'\overline{\otimes} \zeta_{g}':= (\alpha',\beta')}\\
&(X_{f}'' \;\;\;\;\;\;{,}& X_{g}'') &\ar[r]&& X_{f}''\overline{\otimes} X_{g}''}$$

We showed in lemma \ref{zeta is a morphism in both arguments} that $\zeta_{f}\overline{\otimes} \zeta_{g}:= (\alpha,\beta)$ is a morphism in
$MF(K[[x,y]],fg)$, where
 $$(\alpha,\beta)=([\alpha_{f}\otimes \alpha_{g}],[\beta_{f}\otimes \beta_{g}])$$
Similarly, if $\zeta_{f}':=(\alpha_{f}',\beta_{f}')$ and $\zeta_{g}':=(\alpha_{g}',\beta_{g}')$ then $\zeta_{f}' \overline{\otimes} \zeta_{g}'=(\alpha',\beta')$ where $$(\alpha',\beta')=([\alpha_{f}'\otimes \alpha_{g}'],[\beta_{f}'\otimes \beta_{g}'])$$
 It now remains to show the composition and the identity axioms.\\\\
\textit{Identity Axiom}: \\
We show that $F(id_{(X_{f},X_{g})})=id_{F(X_{f},X_{g})}$.\\
Now, $F(id_{(X_{f},X_{g})})=F(id_{X_{f}},id_{X_{g}}):=id_{X_{f}}\overline{\otimes} id_{X_{g}}: X_{f} \overline{\otimes}X_{g} \rightarrow X_{f} \overline{\otimes}X_{g}$.\\
And by definition \ref{defn zeta is a bifuntor}, $id_{X_{f}}\overline{\otimes} id_{X_{g}}$ is the pair of matrices \\
 $([I_{n}\otimes I_{m}],[I_{n}\otimes I_{m}])\,\,\,\,\,\,\,\,\,\,\,\dag$ \\

Next, we compute $id_{F(X_{f},X_{g})}=id_{X_{f} \overline{\otimes}X_{g}}:X_{f} \overline{\otimes}X_{g} \rightarrow X_{f} \overline{\otimes}X_{g}$.\\
By definition of a morphism in the category $MF(fg)$, we know that \\

$id_{X_{f} \overline{\otimes}X_{g}}:=([I_{nm}],[I_{nm}])\,\,\,\,\,\,\,\,\,\,\,\dag\dag
$\\

Since $I_{n}\otimes I_{m}= I_{nm}$, we see that $\dag$ and $\dag\dag$ are the same, therefore $F(id_{(X_{f},X_{g})})=id_{F(X_{f},X_{g})}$ as desired.\\\\
\textit{Composition Axiom}:\\
Consider the situation:
$$\xymatrix@R=4in{X_{f}\ar[r]^{\zeta_{f}} & X_{f}'\ar[r]^{\zeta_{f}'}& X_{f}''}$$ $$\xymatrix@R=4in{X_{g}\ar[r]^{\zeta_{g}} & X_{g}'\ar[r]^{\zeta_{g}'}& X_{g}''}$$
$$\xymatrix@R=4in{X_{f} \overline{\otimes}X_{g}\ar[r]^{F(\zeta_{f},\zeta_{g})} & X_{f}'\overline{\otimes}X_{g}'\ar[r]^{F(\zeta_{f}',\zeta_{g}')}& X_{f}''\overline{\otimes} X_{g}''}$$

We need to show $F(\zeta_{f}'\circ \zeta_{f},\zeta_{g}'\circ \zeta_{g})=F(\zeta_{f}',\zeta_{g}')\circ F(\zeta_{f},\zeta_{g})$. \\
Now, $\zeta_{f}'\circ \zeta_{f}= (\alpha_{f}'\alpha_{f},\beta_{f}'\beta_{f})$ and $\zeta_{g}'\circ \zeta_{g}=(\alpha_{g}'\alpha_{g},\beta_{g}'\beta_{g})$.\\\\
Thanks to definition \ref{defn zeta is a bifuntor}, we obtain:
$$(\zeta_{f}'\circ \zeta_{f})\overline{\otimes}(\zeta_{g}'\circ \zeta_{g})=([\alpha_{f}'\alpha_{f}\otimes \alpha_{g}'\alpha_{g}],[\beta_{f}'\beta_{f}\otimes \beta_{g}'\beta_{g}])\,\,\,\,\,\,\,\,\,\,\,\, \ddag'$$

Next, \\\\
$(\zeta_{f}'\overline{\otimes}\zeta_{g}')\circ (\zeta_{f}\overline{\otimes} \zeta_{g})\\
=([\alpha_{f}'\otimes \alpha_{g}'],[\beta_{f}'\otimes \beta_{g}'])
\circ
([\alpha_{f}\otimes \alpha_{g}],[\beta_{f}\otimes \beta_{g}])
\\
=([\alpha_{f}'\alpha_{f}\otimes \alpha_{g}'\alpha_{g}],[\beta_{f}'\beta_{f}\otimes \beta_{g}'\beta_{g}])
\,\,\,\,\,\,\,\,\,\,\,\, \ddag\ddag'$

From $\ddag'$ and $\ddag\ddag'$, we see that $F(\zeta_{f}'\circ \zeta_{f},\zeta_{g}'\circ \zeta_{g})=F(\zeta_{f}',\zeta_{g}')\circ F(\zeta_{f},\zeta_{g})$.
Thus, $(-)\overline{\otimes}(-)$ is a bifunctor.

\end{enumerate}
\end{proof}

\section{Properties of the refined multiplicative tensor product of matrix factorizations}
In this section, we prove that $\overline{\otimes}$ is associative, commutative and distributive.
\\
We denote by $X_{1}=(\phi_{1},\psi_{1})$ (resp. $X_{2}=(\phi_{2},\psi_{2})$) an $(n_{1}\times n_{1})$ (resp. $(n_{2}\times n_{2})$) matrix factorization of $f\in K[[x]]$. We also let $X'=(\phi',\psi')$ (resp. $X''=(\phi'',\psi'')$) denotes a $(p\times p)$ (resp. $(m\times m)$) matrix factorization of $g\in K[[y]]$ (resp. of $h\in K[[z]]:= K[[z_{1},\cdots,z_{l}]]$).
$X=(\phi,\psi)$ will also be an $r\times r$ matrix factorization of $f\in K[[x]]$.

\subsection{Associativity, commutativity and distributivity of $\overline{\otimes}$ }

\begin{proposition} (Associativity) \label{Assoc of new tensor product}\\
There is an identity:\\
  $(X\overline{\otimes}X')\overline{\otimes}X''=X\overline{\otimes}(X'\overline{\otimes}X'')$ in $MF(fgh)$.

\end{proposition}
\begin{proof}
The desired identity follows from the fact that the standard tensor product for matrices is associative.
\end{proof}
To prove the commutativity of $\overline{\otimes}$, it is worth recalling (cf. section 3.1 \cite{henderson1981vec}) that given two matrices $A$ and $B$, the tensor products $A \otimes B$ and $B \otimes A$ are \textbf{permutation equivalent}. That is, there exist permutation matrices $P$ and $Q$ (so called commutation matrices) such that:
$A \otimes B = P (B \otimes A) Q$. If $A$ and $B$ are square matrices, then $A \otimes B$ and $B\otimes A$ are even \textbf{permutation similar}, meaning we can take $P = Q^{T}$.\\
To be more precise \cite{henderson1981vec}, if $A$ is a $p\times q$ matrix and $B$ is an $r\times s$ matrix, then $$B \otimes A = S_{p,r} (A \otimes B) S_{q,s}^{T}$$
where, $$ S_{m,n}=\sum_{i=1}^{m}( e_{i}^{T}\otimes I_{n}\otimes e_{i}) = \sum_{j=1}^{n}( e_{j}\otimes I_{m}\otimes e_{j}^{T})$$
$I_{n}$ is the $n\times n$ identity matrix and $e_{i}$ is the $i^{th}$ unit vector. $ S_{m,n}$ is the \textbf{perfect shuffle} permutation matrix.
\\
The commutativity of $\overline{\otimes}$ is up to isomorphism. This isomorphism comes from the permutation similarity
of the matrices involved.
\begin{proposition} (commutativity) \label{prop commutativity of new tensor product}\\
For matrix factorizations $X\in MF(f)$ and $X'\in MF(g)$, there is a natural isomorphism

   $X\overline{\otimes} X'\cong X' \overline{\otimes} X \,in\,MF(fg).$

\end{proposition}
\begin{proof}
%
We know that $X\overline{\otimes} X' =([\phi\otimes \phi'],[\psi\otimes \psi'])\,
and\,X'\overline{\otimes} X=([\phi'\otimes \phi],[\psi'\otimes \psi])$.
The desired isomorphism then follows from the fact that $\phi\otimes \phi'$ (respectively $\psi\otimes \psi'$) and $\phi'\otimes \phi$ (respectively $\psi'\otimes \psi$) are permutation similar.
\end{proof}
\begin{proposition} (Distributivity)\\
If $X_{1}$ and $X_{2}$ are matrix factorizations (of $f\in K[[x]]$) of the same size, then there are identities
\begin{enumerate}
  \item $(X_{1}\oplus X_{2})\overline{\otimes} X'=(X_{1}\overline{\otimes} X')\oplus (X_{2}\overline{\otimes} X').$
  \item $ X' \overline{\otimes}(X_{1}\oplus X_{2})=(X'\overline{\otimes}X_{1})\oplus (X'\overline{\otimes}X_{2}).$
  \end{enumerate}

\end{proposition}
\begin{proof}
\begin{enumerate}
  \item
$$(X_{1}\overline{\otimes} X')\oplus (X_{2}\overline{\otimes} X')$$
$\,\,\,\,\,\,\,\,\,\,\,\,\,\,\,\,\,\,\,\,\,\,\,\,\,\,\,\,\,\,\,\,\,\,\,\,\,\,=([\phi_{1}\otimes \phi'],[\psi_{1}\otimes \psi'])\oplus ([\phi_{2}\otimes \phi'],[\psi_{2}\otimes \psi'])$\\
\[=(\begin{bmatrix}
  \phi_{1}\otimes \phi' &   0      \\
  0   & \phi_{2}\otimes \phi'
\end{bmatrix}, \begin{bmatrix}
  \psi_{1}\otimes \psi' &   0     \\
  0   & \psi_{2}\otimes \psi'
\end{bmatrix})\,\,\,\,\,\,\,\,\cdots (\sharp)\]\\
Next, $$(X_{1}\oplus X_{2})\overline{\otimes} X'=((\phi_{1},\psi_{1})\oplus (\phi_{2},\psi_{2}))\overline{\otimes} (\phi',\psi')$$
\[=(\begin{bmatrix}
  \phi_{1} &          0      \\
    0               &   \phi_{2}
\end{bmatrix}, \begin{bmatrix}
  \psi_{1}   &          0      \\
    0               &   \psi_{2}
\end{bmatrix})\overline{\otimes}(\phi',\psi')\]\\
\[=(
  \begin{bmatrix}
  \phi_{1} &          0      \\
    0               &   \phi_{2}
\end{bmatrix}\otimes\phi' , \begin{bmatrix}

  \psi_{1}   &          0      \\
    0               &   \psi_{2}
\end{bmatrix}\otimes\psi')\]\\
\[=(\begin{bmatrix}
  \phi_{1}\otimes \phi' &   0      \\
  0   & \phi_{2}\otimes \phi'
\end{bmatrix}, \begin{bmatrix}
  \psi_{1}\otimes \psi' &   0      \\
  0   & \psi_{2}\otimes \psi'
\end{bmatrix})\,\,\,\,\,\,\,\,\cdots (\sharp')\]
\\%
The desired identity now follows from $(\sharp)$ and $(\sharp')$.
\item The proof of this equality is similar to the foregoing proof.
\end{enumerate}
\end{proof}

In the next section, we use the \textit{refined multiplicative tensor product of matrix factorizations} $\overline{\otimes}$ to refine the algorithm for factoring polynomials proposed in \cite{fomatati2022tensor} on the class  of \textit{summand-reducible polynomials} (cf. definition \ref{defn summand reducible polynomials}).

\section{A refined algorithm for matrix factorization of polynomials on the class of summand-reducible polynomials} \label{chap: improvement of the std method}

We first recall a standard algorithm for factoring polynomials which dates to the 1980s when Kn\"{o}rrer exploited it to prove his celebrated periodicity theorem (cf. theorem 2.1 \cite{brown2016knorrer}). This standard technique, usually referred to as the standard method \cite{crisler2016matrix} for factoring polynomials, builds matrix factorizations of sums of polynomials from "factorizations" of their summands. One conspicuous downside of this algorithm is that for each new summand that is added to the polynomial being factorized, the size (i.e., the number of rows and columns) of the matrix factorization doubles. \\
In \cite{fomatati2022tensor}, the standard method for factoring polynomials was improved on the class of \textit{summand-reducible polynomials} in the sense that matrix factors yielded by this improved algorithm are smaller in size. The main ingredients in that improved algorithm are the multiplicative tensor product of matrix factorization $\widetilde{\otimes}$ and the Yoshino tensor product of matrix factorizations $\widehat{\otimes}$. In this section, we show that if we replace $\widetilde{\otimes}$ by its refined version $\overline{\otimes}$ in that improved algorithm, we obtain a refined algorithm which yields matrix factors that are smaller in size as compared to the ones obtained with the improved algorithm of \cite{fomatati2022tensor}.\\
 In fact, a \textit{summand-reducible polynomial} is one that can be written in the form $f=t_{1}+\cdots + t_{s}+ g_{11}\cdots g_{1m_{1}} + \cdots + g_{l1}\cdots g_{lm_{l}}$ under some specified conditions where each $t_{k}$ is a monomial and each $g_{ji}$ is a sum of monomials. In \cite{fomatati2022tensor}, it is proved that if $p_{ji}$ is the number of monomials in $g_{ji}$, then there is an improved version of the standard method for factoring $f$ which produces factorizations of size $2^{\prod_{i=1}^{m_{1}}p_{1i} + \cdots + \prod_{i=1}^{m_{l}}p_{li} - (\sum_{i=1}^{m_{1}}p_{1i} + \cdots + \sum_{i=1}^{m_{l}}p_{li})}$ times smaller than the size one would normally obtain with the standard method. In this paper, we show that under the same hypothesis, the refined algorithm that we are going to present produces factorizations of size $2^{\prod_{i=1}^{m_{1}}p_{1i} + \cdots + \prod_{i=1}^{m_{l}}p_{li} - (\sum_{i=1}^{m_{1}}p_{1i} + \cdots + \sum_{i=1}^{m_{l}}p_{li})+\sum_{j=1}^{l}m_{j}-l}$ times smaller than the size one would normally obtain with the standard method. As we will show below, this means that the matrix factors we obtain with the refined algorithm are $2^{\sum_{j=1}^{l}m_{j}-l}$ times smaller than the ones obtained with the improved method presented in \cite{fomatati2022tensor}.
\\
In our presentation, we limit ourselves to polynomials in the ring $S=\mathbb{R}[x_{1}\dots,x_{n}]$ where $\mathbb{R}$ is the set of real numbers.

\subsection{The standard method for factoring polynomials} \label{sec: the std method}
\textbf{Introduction}\\
In his seminal paper \cite{eisenbud1980homological}, Eisenbud proved that using matrices, both reducible and irreducible polynomials in S can be factorized. He showed that the matrix factorizations of the polynomial $f$ are intimately related
to homological properties of modules over the quotient ring $S/(f)$, known as the hypersurface ring. \cite{knorrer1987cohen} and \cite{buchweitz1987cohen} contain more background on the connection between matrix factorizations and algebraic geometry. These papers have details on the connection that exists between matrix factorizations and maximal Cohen-Macaulay Modules. In this subsection, we describe a way to construct matrix factorizations of a polynomial without resorting to the homological methods that Eisenbud introduced.\\\\
\textbf{The standard method}\\
Here, we recall the standard technique for factoring polynomials using matrices.
\begin{proposition} \cite{crisler2016matrix}
  For $i,j\in \{1,2\}$, let $(C_{i},D_{i})$ denote an $n\times n$ matrix factorization of the polynomial $f_{i}\in S$. In addition, assume that the matrices $C_{i}$ and $D_{j}$ commute when $i\neq j.$ Then the matrices

  $$
\begin{pmatrix}
  C_{1}
  & \rvline & -D_{2} \\
\hline
  C_{2} & \rvline &
  D_{1}
\end{pmatrix},
\begin{pmatrix}
  D_{1}
  & \rvline & D_{2} \\
\hline
  -C_{2} & \rvline &
  C_{1}
\end{pmatrix}
$$
give a $2n\times 2n$ matrix factorization of $f_{1} + f_{2}$.
\end{proposition}
 The following consequence of the foregoing result is actually the basis for the main construction of the standard algorithm for matrix factorization of polynomials.
 \begin{corollary}\cite{crisler2016matrix} \label{cor: standard method for factoring polyn}
 If $(C,D)$ is an $n\times n$ matrix factorization of $f$ and $g,h$ are two polynomials, then
 $$
\begin{pmatrix}
  C
  & \rvline & -gI_{n} \\
\hline
  hI_{n} & \rvline &
  D
\end{pmatrix},
\begin{pmatrix}
  D
  & \rvline & gI_{n} \\
\hline
  -hI_{n} & \rvline &
  C
\end{pmatrix}
$$

give a $2n\times 2n$ matrix factorization of $f + gh$.
 \end{corollary}
 \begin{proof}
   Since the matrices $gI_{n}$ and $hI_{n}$ commute with all $n\times n$ matrices, the proof follows from the previous proposition.
 \end{proof}

 Thanks to this corollary, one can inductively construct matrix factorizations of polynomials of the form:\\
 $$f = f_{k} = g_{1}h_{1} + g_{2}h_{2}+ \cdots + g_{k}h_{k}.$$
 For $k=1$, we have $f=g_{1}h_{1}$ and clearly $[g_{1}][h_{1}]=[g_{1}h_{1}]=[f_{1}]$ is a $1\times 1$ matrix factorization. Next, assume that $C$ and $D$ are matrix factorizations of $f_{k-1}$, i.e., $CD=If_{k-1}$ where $I$ is the identity matrix of the right size. Hence, using the foregoing corollary, we obtain a matrix factorization of $f_{k}$:
 $$
(\begin{pmatrix}
  C
  & \rvline & -g_{k}I_{n} \\
\hline
  h_{k}I_{n} & \rvline &
  D
\end{pmatrix},
\begin{pmatrix}
  D
  & \rvline & g_{k}I_{n} \\
\hline
  -h_{k}I_{n} & \rvline &
  C
\end{pmatrix})
$$
 \begin{definition} \cite{crisler2016matrix}
   The foregoing algorithm is called \textbf{the standard method} for factoring polynomials.
 \end{definition}
 Other synonyms for standard method are \textit{standard technique} or \textit{standard algorithm}.
 \begin{remark} 
Let $(P,Q)$ be the matrices in Corollary \ref{cor: standard method for factoring polyn}, we observe that if we interchange the rows (respectively the columns) of $P$ and interchange the columns (respectively the rows) of $Q$, then we would still obtain a $2n\times 2n$ matrix factorization of $f + gh$.

 \end{remark}

 Since every polynomial can be expressed as a sum of finitely many monomials, the standard method can be used to produce matrix factorizations of any polynomial.\\
 Though this algorithm works for any polynomial, it has a noticeable downside. The sizes of factorizations grow very quickly due to the fact that for every new summand $g_{n}h_{n}$ added to the polynomial, the factorizations double in size. It is easy to see that with this method, to factor a polynomial with $k$ summands, say $$f_{k}=g_{1}h_{1} + g_{2}h_{2}+ \cdots + g_{k}h_{k},$$ one obtains matrices of size $2^{k-1}$. Thus, the size of matrix factors can grow extremely large very quickly. For example if $k=8$, we will obtain matrices of size $ 2^{7}= 128$ and for $k=10$, we will obtain matrices of size $2^{9} = 512$.
\\The following two examples that illustrate the standard method show that if a polynomial is written in its expanded form (sum of monomials), then the matrix factors have monomial entries. But if it is not in its expanded form, then some entries of the matrix factors will not be monomials but sums of monomials.
%
 \begin{example} \label{exple: good matrix facto of g}
   Let $g= xy^{2}+x^{2}z+yz^{2}$. We use the standard method to find a matrix factorization of $g$.
   First a matrix factorization of $xy^{2}+x^{2}z$ is
    \[
  (\begin{bmatrix}
    x & -x^{2} \\
    z & y^{2}
  \end{bmatrix},\begin{bmatrix}
    y^{2} & x^{2} \\
    -z & x
  \end{bmatrix})
\]
Hence, a matrix factorization of $g= xy^{2}+x^{2}z+yz^{2}$ is then:
 \[
  P=(\begin{bmatrix}
    x & -x^{2} & -y & 0\\
    z &   y^{2}    &  0 & -y\\
 z^{2} &  0    &  y^{2} &  x^{2}\\
     0 &  z^{2}&  -z & x
  \end{bmatrix},\begin{bmatrix}
     y^{2} & x^{2} &  y & 0\\
    -z &   x   &  0 & y\\
 -z^{2} &  0    &  x &  -x^{2}\\
     0 &  -z^{2}&  z & y^{2}
  \end{bmatrix})
\]
\end{example}
 \begin{example} \label{exple: bad matrix facto of g}
 Let $l= xy^{2}+(x^{2}+yz)z$. Observe that $l=g$ where $g$ is given in example \ref{exple: good matrix facto of g}. We use the standard method and quickly find a matrix factorization of $l$:

  \[
  Q=(\begin{bmatrix}
    x & -(x^{2}+yz) \\
    z & y^{2}
  \end{bmatrix},\begin{bmatrix}
    y^{2} & (x^{2}+yz) \\
    -z & x
  \end{bmatrix})
\]
 \end{example}
We observe that the factorization we obtain for $l$ is not as nice as the one we obtain for $g$, in the sense that the complexity of some entries in the factorization of $l$ is higher than what we have for $g$.
For instance, in $Q$ the entry $(x^{2}+yz)$ is more complex than all the entries in $P$ of example \ref{exple: good matrix facto of g}. This shows that it is better to use the expanded version of a polynomial to find its matrix factorization.\\
%
%
 \begin{quote}
   So, just like in \cite{fomatati2022tensor} we make the important assumption that before applying the standard method to a given polynomial, it has to be written in its expanded form.
 \end{quote}
     But this comes at a price! The size of the factorizations becomes big as one could notice in  examples \ref{exple: good matrix facto of g} and \ref{exple: bad matrix facto of g}\\
In \cite{fomatati2022tensor}, it was proved that by combining the Yoshino tensor product $\widehat{\otimes}$ and the multiplicative tensor product $\widetilde{\otimes}$, one can obtain smaller matrix factors than the ones produced by the standard method.\\
 In the following section, we will show that the refined multiplicative tensor product of matrix factorizations $\overline{\otimes}$ can be used in place of the multiplicative tensor product of matrix factorizations $\widetilde{\otimes}$ to produce even better results on the size of the matrix factors.

\subsection{The refined algorithm} \label{sec: the refined algorithm}
In this section, we use the refined multiplicative tensor product of matrix factorizations $\overline{\otimes}$ in combination with the tensor product of matrix factorizations $\widehat{\otimes}$ to refine the improved algorithm (see section 4 of \cite{fomatati2022tensor}) for matrix factorization of polynomials on the class of \textit{summand-reducible} polynomials (cf. Definition \ref{defn summand reducible polynomials}). In fact, we show that our approach produces factorizations that are of smaller sizes than the factorizations produced by the improved algorithm in \cite{fomatati2022tensor} on the aforementioned class of polynomials.\\

In this paper, we will refer to the algorithm developed in the proof of Theorem \ref{thm improved algo for summand-red polyn}  (cf. Theorem 4.2 of \cite{fomatati2022tensor}) as
 \textbf{the improved algorithm} and we will refer to the one developed in the proof of Theorem \ref{thm refined algo for summand-red polyn} as
\textbf{the refined algorithm}.\\

We now recall the definition of the class of \textit{summand-reducible} polynomials which is made up of polynomials in which some monomials can be factorized in a nice way, hence allowing the polynomial to be written with less summands.
\begin{definition}\cite{fomatati2022tensor} \label{defn summand reducible polynomials}
  A polynomial $f$ is said to be \textbf{summand-reducible} if it can be written in the form $$f=t_{1}+\cdots + t_{s}+ g_{11}\cdots g_{1m_{1}} + \cdots + g_{l1}\cdots g_{lm_{l}},$$ where:
  \begin{enumerate}
    \item $\bullet$ If $s=0$, then there exist at least two products $g_{11}\cdots g_{1m_{1}}$ and $g_{21}\cdots g_{2m_{2}}$ in $f$.\\
        $\bullet$ If $s\neq 0$, then there exists at least one product $g_{11}\cdots g_{1m_{1}}$ in $f$.
    \item For $i=1,\cdots,s$; each $t_{i}$ is a monomial and so $t_{i}=h_{i1}h_{i2}$, where $h_{i1}$ and $h_{i2}$ are products of variables possibly raised to some power.
    \item For $j=1,\cdots, l$; each $g_{j1}\cdots g_{jm_{j}}$ is a product of sums of monomials, such that if it is expanded, $g_{j1}\cdots g_{jm_{j}}$ would have more monomials than the number that appears in the factor form $g_{j1}\cdots g_{jm_{j}}$.
    \item For $1\leq j \leq l$, at least one of the products $g_{j1}\cdots g_{jm_{j}}$ has at least two factors.
  \end{enumerate}
\end{definition}
\begin{example} \label{exples of simple summand-red polyn}
  The following are summand-reducible polynomials:
  \begin{enumerate}
    \item $z^{5}+ yx^{3}+ zx^{4} + yz^{2}x^{2} + xy^{3} + x^{2}zy^{2} + z^{2}y^{3}= z^{5}+ (xy + x^{2}z + yz^{2})(x^{2} + y^{2})$.
    \item  $x^{5}-y^{5}+ (xy + yz^{2})(x^{2} + y^{2} + z)$.
    \item $z^{3} + xyz + yx^{3}+ zx^{4} + yz^{2}x^{2} + xy^{3} + x^{2}zy^{2} + z^{2}y^{3} + xyz + x^{2}z^{2} + yz^{3}=z^{3} + xyz + (xy + x^{2}z + yz^{2})(x^{2} + y^{2} + z)$.
  \end{enumerate}
\end{example}
\begin{remark}
  Observe that our definition \ref{defn summand reducible polynomials} mostly targets polynomials with more than six monomials because factorizations obtained with the standard method begin to be of considerable sizes.
\end{remark}
Consider the following polynomials for which some conditions of definition \ref{defn summand reducible polynomials} fail:
\begin{enumerate}
  \item Let $n,m\in \mathbb{N}$, $x^{m}-y^{n}$ is not summand-reducible. Here, the first condition fails.
  \item $zx+(x-y)(x^{4} + x^{3}y + x^{2}y^{2}+xy^{3}+y^{4})= zx+x^{5}-y^{5}$ is not summand-reducible. Here, the third condition fails. Note that it is better to write $(x-y)(x^{4} + x^{3}y + x^{2}y^{2}+xy^{3}+y^{4})$ as $x^{5}-y^{5}$ because the latter expression produces smaller factorizations than the former (this can be verified when reading the proof of theorem \ref{thm refined algo for summand-red polyn}).
\end{enumerate}

\begin{definition} \cite{fomatati2022tensor}
  A polynomial $f$ is said to be \textbf{summand-reduced} if it is in the form $f=t_{1}+\cdots + t_{s}+ g_{11}\cdots g_{1m_{1}} + \cdots + g_{l1}\cdots g_{lm_{l}}$ described in definition \ref{defn summand reducible polynomials}.
\end{definition}
\begin{example}
   $f= zy +(xy^{2}+x^{2}z+yz^{2})(xy+z^{2}) +(yz + xy^{2} +x^{2})(x^{3}z^{2}+ yx + y^{2})$ is a summand-reduced polynomial.
  \end{example}

\begin{remark}
  With the standard method, even if one knows matrix factorizations of polynomials $f$ and $g$, one cannot derive from them a matrix factorization of the product $fg$ nor of the sum $f+g$.\\
One of the main ingredients used in the improved algorithm (section 4 of \cite{fomatati2022tensor}) is Yoshino's tensor product of matrix factorizations $\widehat{\otimes}$, because it produces a matrix factorization of the sum of two polynomials from the matrix factorizations of each of these polynomials.\\
Another crucial ingredient used in the improved algorithm is the multiplicative tensor product of matrix factorizations $\widetilde{\otimes}$ (recalled in definition \ref{defn of the multiplicative tensor product}) which produces a matrix factorization of the product of two polynomials from the matrix factorizations of each of these polynomials.\\
The proofs of theorems 4.1 and 4.2 of \cite{fomatati2022tensor} show how these two bifunctorial operations ($\widehat{\otimes}$ and $\widetilde{\otimes}$) help in reducing the size of matrix factors of summand-reducible polynomials.
\end{remark}
In the sequel, we will recall Theorem 4.2 of \cite{fomatati2022tensor} which is the main result (i.e., the improved algorithm result) we aim at further improving in this paper. Next, we will somehow rewrite the improved algorithm to obtain a refined algorithm which will help obtain matrix factors that are smaller in size as compared to the ones obtained by the improved algorithm. The main operation we will perform in the improved algorithm in order to obtain a refined one, will be to replace the multiplicative tensor product of matrix factorizations $\widetilde{\otimes}$ by the refined multiplicative tensor product of matrix factorizations $\overline{\otimes}$.
\\
We recall the following result which will be refined in Theorem \ref{thm refined algo for summand-red polyn} :\\
\begin{theorem} \cite{fomatati2022tensor} \label{thm improved algo for summand-red polyn}
Let $f=t_{1}+\cdots + t_{s}+ g_{11}\cdots g_{1m_{1}} + \cdots + g_{l1}\cdots g_{lm_{l}}$ be a summand-reducible polynomial. Let $p_{ji}$ be the number of monomials in $g_{ji}$. Then
there is an improved version of the standard method for factoring $f$ which produces factorizations of size $$2^{\prod_{i=1}^{m_{1}}p_{1i} + \cdots + \prod_{i=1}^{m_{l}}p_{li} - (\sum_{i=1}^{m_{1}}p_{1i} + \cdots + \sum_{i=1}^{m_{l}}p_{li})}$$ times smaller than the size one would normally obtain with the standard method.
\end{theorem}
In the proof of this theorem (cf. proof of Theorem 4.2 \cite{fomatati2022tensor}), it was shown that for a summand-reduced polynomial $f$ as defined in the theorem, the improved algorithm produces matrix factors of sizes $2^{\sum_{i=1}^{m_{1}}p_{1i} + \cdots + \sum_{i=1}^{m_{l}}p_{li} +s-1}$.\\
\textbf{Nota Bene:} The proof of Theorem \ref{thm refined algo for summand-red polyn} shows how to combine the operations $\widehat{\otimes}$ and $\overline{\otimes}$ to obtain matrix factors for a given summand-reduced polynomial $f$,  whereas the proof of Theorem \ref{thm improved algo for summand-red polyn}
shows how to combine the operations
$\widehat{\otimes}$ and $\widetilde{\otimes}$ to obtain matrix factors of $f$. Though $\widetilde{\otimes}$ is simply the direct sum of two copies of $\overline{\otimes}$, the matrix factors one obtains when combining $\widehat{\otimes}$ and $\overline{\otimes}$ are not simply as twice as small as the ones obtained when combining $\widehat{\otimes}$ and $\widetilde{\otimes}$. In fact, it depends on the number of terms $f$ has. This is shown in Corollary \ref{Cor: comparison between the two algorithms}.\\

We can now state and prove the following theorem:
\begin{theorem} \label{thm refined algo for summand-red polyn}
Let $f=t_{1}+\cdots + t_{s}+ g_{11}\cdots g_{1m_{1}} + \cdots + g_{l1}\cdots g_{lm_{l}}$ be a summand-reducible polynomial. Let $p_{ji}$ be the number of monomials in $g_{ji}$. Then
there is an improved version of the standard method for factoring $f$ which produces factorizations of size $$2^{\prod_{i=1}^{m_{1}}p_{1i} + \cdots + \prod_{i=1}^{m_{l}}p_{li} - (\sum_{i=1}^{m_{1}}p_{1i} + \cdots + \sum_{i=1}^{m_{l}}p_{li} -  \sum_{j=1}^{l}m_{j} + l)}$$ times smaller than the size one would normally obtain with the standard method.
\end{theorem}

\begin{proof}
  First, we construct the algorithm, then we prove that the resulting matrix factorizations (for summand-reducible polynomials) are$$2^{\prod_{i=1}^{m_{1}}p_{1i} + \cdots + \prod_{i=1}^{m_{l}}p_{li} - (\sum_{i=1}^{m_{1}}p_{1i} + \cdots + \sum_{i=1}^{m_{l}}p_{li} -  \sum_{j=1}^{l}m_{j} + l)}$$ times smaller in size than what one would obtain with the standard method.\\
  We inductively construct the matrix factorizations of summand-reduced polynomials using the tensor products $\widehat{\otimes}$ and $\overline{\otimes}$  that were not existing in the 1980s when the standard method was developed.
  The algorithm we propose here is a refinement of the one given in \cite{fomatati2022tensor} for summand-reducible polynomials (cf. proof of theorem 4.2 of \cite{fomatati2022tensor}).

Let $f=t_{1}+\cdots + t_{s}+ g_{11}\cdots g_{1m_{1}} + \cdots + g_{l1}\cdots g_{lm_{l}}$ be a summand-reducible polynomial. Let $p_{ji}$ be the number of monomials in $g_{ji}$. \\
If $\forall k\in \{1,\cdots,s\}$, $t_{k}= 0$, then do:
\begin{enumerate}
  \item For each $j\in \{1,\cdots, l\}$ and $i\in \{1,\cdots, m_{j}\}$, use the standard method to find a matrix factorization of $g_{ji}$ of size $2^{p_{ji}-1}$.
  \item Next, for each $j\in \{1,\cdots, l\}$; use the refined multiplicative tensor product of matrix factorizations $\overline{\otimes}$ to find a matrix factorization of $g_{j1}\cdots g_{jm_{j}}$ of size $$2^{\sum_{i=1}^{m_{j}}p_{ji}-m_{j}}$$
  \item Now use the tensor product of matrix factorizations $\widehat{\otimes}$ to find a matrix factorization of $g_{11}\cdots g_{1m_{1}} + \cdots + g_{l1}\cdots g_{lm_{l}}$ of size
      $$(2^{l-1})(\prod_{j=1}^{l}(2^{\sum_{i=1}^{m_{j}}p_{ji}-m_{j}})= 2^{l-1 + \sum_{i=1}^{m_{1}}p_{1i} + \cdots + \sum_{i=1}^{m_{l}}p_{li} - \sum_{j=1}^{l}m_{j}}.$$
      Let us find the size of matrix factors the standard method would produce for $$g_{11}\cdots g_{1m_{1}} + \cdots + g_{l1}\cdots g_{lm_{l}}.$$
      Let $n_{j}=$number of monomials in the expanded form of the $j^{th}$ product $g_{j1}\cdots g_{jm_{j}}$. Then $n_{j}=\prod_{i=1}^{m_{j}} p_{ji}$. Hence, the number of monomials in the expanded form of $g_{11}\cdots g_{1m_{1}} + \cdots + g_{l1}\cdots g_{lm_{l}}$ would be $\sum_{j=1}^{l} n_{j}= \sum_{j=1}^{l} \prod_{i=1}^{m_{j}} p_{ji}$. \\So the size of factorizations produced by the standard method would be $2^{(\sum_{j=1}^{l} \prod_{i=1}^{m_{j}}p_{ji})-1}.$ \\
      Thus, the size of matrix factors produced by our refined algorithm would be
       $$2^{(\sum_{j=1}^{l} \prod_{i=1}^{m_{j}}p_{ji})-1}\div 2^{l-1 + \sum_{i=1}^{m_{1}}p_{1i} + \cdots + \sum_{i=1}^{m_{l}}p_{li} - \sum_{j=1}^{l}m_{j}}= 2^{(\sum_{j=1}^{l} \prod_{i=1}^{m_{j}}p_{ji})- (l + \sum_{i=1}^{m_{1}}p_{1i} + \cdots + \sum_{i=1}^{m_{l}}p_{li} - \sum_{j=1}^{l}m_{j})}$$
      times smaller in size than the factorizations produced by the standard method.
  \item If there exists $k\in \{1,\cdots,s\}$ such that $t_{k}\neq 0$, then use the standard method to inductively find a matrix factorization $(A,B)$ of $t_{1}+\cdots + t_{s}$ of size $2^{s-1}$.
  \item Then do steps 1), 2) and 3) above to find a matrix factorization $(C,D)$ of $g_{11}\cdots g_{1m_{1}} + \cdots + g_{l1}\cdots g_{lm_{l}}$ of size $ 2^{l-1 + \sum_{i=1}^{m_{1}}p_{1i} + \cdots + \sum_{i=1}^{m_{l}}p_{li} - \sum_{j=1}^{l}m_{j}}$.
  \item Now, use $\widehat{\otimes}$ to find a matrix factorization $(A,B) \widehat{\otimes} (C,D)$ of $f=t_{1}+\cdots + t_{s}+ g_{11}\cdots g_{1m_{1}} + \cdots + g_{l1}\cdots g_{lm_{l}}$ of size $$2(2^{s-1})( 2^{l-1 + \sum_{i=1}^{m_{1}}p_{1i} + \cdots + \sum_{i=1}^{m_{l}}p_{li} - \sum_{j=1}^{l}m_{j}})= 2^{l-1 + \sum_{i=1}^{m_{1}}p_{1i} + \cdots + \sum_{i=1}^{m_{l}}p_{li} - \sum_{j=1}^{l}m_{j}+s}.$$
      Note that $f$ in expanded form has $$\sum_{j=1}^{l} n_{j}+ s= (\sum_{j=1}^{l} \prod_{i=1}^{m_{j}}p_{ji}) + s$$ monomials and so the standard method would produce factorizations of size\\ $2^{(\sum_{j=1}^{l} \prod_{i=1}^{m_{j}}p_{ji}) + s - 1}$.\\ Hence the factorizations our refined algorithm produces are $$2^{(\sum_{j=1}^{l} \prod_{i=1}^{m_{j}}p_{ji}) + s - 1}\div 2^{l-1 + \sum_{i=1}^{m_{1}}p_{1i} + \cdots + \sum_{i=1}^{m_{l}}p_{li} - \sum_{j=1}^{l}m_{j}+s}= 2^{(\sum_{j=1}^{l} \prod_{i=1}^{m_{j}}p_{ji})- (\sum_{i=1}^{m_{1}}p_{1i} + \cdots + \sum_{i=1}^{m_{l}}p_{li})+\sum_{j=1}^{l}m_{j}-l}$$
      times smaller in size than the factorizations produced by the standard method. QED.

\end{enumerate}
\end{proof}
The following corollary is a comparison between the sizes of matrix factors (of a given summand-reducible polynomial) produced by the improved algorithm on the one hand and those produced by the refined algorithm on the other hand.
\begin{corollary} \label{Cor: comparison between the two algorithms}
  Let $f=t_{1}+\cdots + t_{s}+ g_{11}\cdots g_{1m_{1}} + \cdots + g_{l1}\cdots g_{lm_{l}}$ be a summand-reduced polynomial. Let $p_{ji}$ be the number of monomials in $g_{ji}$. Then the refined algorithm produces matrix factors of $f$ whose size is
  $$2^{\sum_{j=1}^{l}m_{j}-l} $$
  times smaller than the size one would normally obtain with the improved algorithm.
\end{corollary}
\begin{proof}
  We know from the proofs of Theorems \ref{thm improved algo for summand-red polyn} and \ref{thm refined algo for summand-red polyn} that the improved and the refined algorithm produce respectively matrix factors of sizes $2^{\sum_{i=1}^{m_{1}}p_{1i} + \cdots + \sum_{i=1}^{m_{l}}p_{li} +s-1}$
 and $2^{l-1 + \sum_{i=1}^{m_{1}}p_{1i} + \cdots + \sum_{i=1}^{m_{l}}p_{li} - \sum_{j=1}^{l}m_{j}+s}$. Hence, the refined algorithm produces matrix factors of $f$ whose size is
  $$2^{\sum_{i=1}^{m_{1}}p_{1i} + \cdots + \sum_{i=1}^{m_{l}}p_{li} +s-1}\div 2^{l-1 + \sum_{i=1}^{m_{1}}p_{1i} + \cdots + \sum_{i=1}^{m_{l}}p_{li} - \sum_{j=1}^{l}m_{j}+s} =2^{\sum_{j=1}^{l}m_{j}-l} $$
  times smaller than the size one would normally obtain with the improved algorithm.
\end{proof}

\begin{example} \textbf{Illustration of Theorem \ref{thm refined algo for summand-red polyn} and Corollary \ref{Cor: comparison between the two algorithms}} \label{exple: first one after the main thm} \\
  Let $f= zy +(xy^{2}+x^{2}z+yz^{2})(xy+z^{2}) +(yz + xy^{2} +x^{2})(x^{3}z^{2}+ yx + y^{2})$.\\
  $f$ in expanded form has $1 + 3\times 2 + 3\times 3=16$ monomials and so the standard method will produce factorizations of size $2^{16-1}=2^{15}$.\\
  From the proof of Theorem \ref{thm refined algo for summand-red polyn}, we know that for a polynomial $$f=t_{1}+\cdots + t_{s}+ g_{11}\cdots g_{1m_{1}} + \cdots + g_{l1}\cdots g_{lm_{l}}$$ the refined algorithm produces matrix factors of size $$2^{l-1 + \sum_{i=1}^{m_{1}}p_{1i} + \cdots + \sum_{i=1}^{m_{l}}p_{li} - \sum_{j=1}^{l}m_{j}+s}.$$
       For this example: $s=1$, $l=2$, $m_{1}=2$, $m_{2}=2$, $p_{11}=p_{21}=p_{22}=3$ and $p_{12}=2$. So, our algorithm would produce factorizations of size $$2^{l-1+p_{11}+p_{12}+p_{21}+p_{22}-m_{1}-m_{2}+s}=2^{2-1+3+2+3+3-2-2+1}=2^{9}.$$
   Hence, from Theorem \ref{thm refined algo for summand-red polyn} we deduce that the refined algorithm produces factorizations of size $2^{15} \div 2^{9}=2^{6}=64$ times smaller than what the standard method produces!
\\
From the proof of Theorem \ref{thm improved algo for summand-red polyn} (cf. Theorem 4.2 of \cite{fomatati2022tensor}), we know that for a summand-reducible polynomial $f=t_{1}+\cdots + t_{s}+ g_{11}\cdots g_{1m_{1}} + \cdots + g_{l1}\cdots g_{lm_{l}}$ the improved algorithm produces matrix factors of size $$ 2^{\sum_{i=1}^{m_{1}}p_{1i} + \cdots + \sum_{i=1}^{m_{l}}p_{li}+s-1}.$$
So, the improved algorithm would produce factorizations of size $$2^{p_{11}+p_{12}+p_{21}+p_{22}+s-1}=2^{3+2+3+3+1-1}=2^{11}.$$\\
  Hence, from Theorem \ref{thm improved algo for summand-red polyn} we deduce that the improved algorithm produces factorizations of size $2^{15} \div 2^{11}=2^{4}=16$ times smaller than what the standard method produces!\\
  We can now illustrate the result of Corollary \ref{Cor: comparison between the two algorithms}, namely that for a summand-reduced polynomial $$f=t_{1}+\cdots + t_{s}+ g_{11}\cdots g_{1m_{1}} + \cdots + g_{l1}\cdots g_{lm_{l}}$$ the refined algorithm produces matrix factors whose size is
  $$2^{\sum_{j=1}^{l}m_{j}-l} $$
  times smaller than the size one would normally obtain with the improved algorithm. In fact,
  for this example, $l=2$, $m_{1}=2$, $m_{2}=2$, thus the refined algorithm produces matrix factors of $f$ whose size is $$2^{\sum_{j=1}^{l}m_{j}-l}=2^{m_{1}+m_{2}-2}=2^{2+2-2}=4$$ times smaller than the size one obtains with the improved algorithm. This confirms what we obtained in the previous paragraphs (in this example) since $4=2^{11}\div 2^{9}$.
\end{example}

\begin{example} \label{exple: matrix factors presented, part II}
Use the refined algorithm to factorize the polynomial $f= x^{5}y^{2} + (xy^{2}+x^{2}z+yz^{2})(x^{2}z+y^{2}+y^{2}z)$ and compare the size of the matrix factors with the one obtained using the standard method.\\
Since $f$ in its expanded form has $1+3\times 3=1+9=10$ monomials, the size of matrix factors obtained using the standard method would be $2^{10-1}=2^{9}=512$. We can use Theorem \ref{thm refined algo for summand-red polyn} as we did in Example \ref{exple: first one after the main thm} to find that the size of matrix factors of $f$ using the refined algorithm is $2^{1-1+3+3-2+1}=2^{5}=32$, that is $\frac{512}{32}=16$ times smaller than the size obtained using the standard method.\\
In the sequel, we are going to use the refined algorithm to find matrix factors of $f$ and we will see that they are actually of size $32$.\\
 Let $g=xy^{2}+x^{2}z+yz^{2}$ and $t=x^{2}z+y^{2}+y^{2}z$, so that $f=x^{5}y^{2} + gt$.\\

 In Example \ref{exple: good matrix facto of g}, we used the standard method to find a matrix factorization of the polynomial $g$:

\[
  P=(\phi_{g},\psi_{g})=(\begin{bmatrix}
    x & -x^{2} & -y & 0\\
    z &   y^{2}    &  0 & -y\\
 z^{2} &  0    &  y^{2} &  x^{2}\\
     0 &  z^{2}&  -z & x
  \end{bmatrix},\begin{bmatrix}
     y^{2} & x^{2} &  y & 0\\
    -z &   x   &  0 & y\\
 -z^{2} &  0    &  x &  -x^{2}\\
     0 &  -z^{2}&  z & y^{2}
  \end{bmatrix})
\]

Let us find a matrix factorization of $t=x^{2}z+y^{2}+y^{2}z=d+y^{2}z$, where $d=x^{2}z+y^{2}$. Using the standard method, we find that
\[
  (\begin{bmatrix}
    x^{2} & -y \\
    y & z
  \end{bmatrix},\begin{bmatrix}
    z & y \\
    -y & x^{2}
  \end{bmatrix})
\]
  is a matrix factorization of $d=x^{2}z+y^{2}$. Thus, using the standard method, a matrix factorization of the polynomial $t$ is

\[
  N=(\phi_{t},\psi_{t})=(\begin{bmatrix}
    x^{2} & -y & -y^{2} & 0\\
    y &  z  &  0 & -y^{2}\\
   z &  0  &  z &  y\\
     0 &  z&  -y & x^{2}
  \end{bmatrix},\begin{bmatrix}
     z & y &  y^{2} & 0\\
    -y &  x^{2}   &  0 & y^{2}\\
 -z &  0    &  x^{2} &  -y\\
     0 &  -z &  y & z
  \end{bmatrix})
\]

According to the proof of Theorem \ref{thm refined algo for summand-red polyn} to find a matrix factorization for $f$, we need to:
\begin{enumerate}
  \item First of all find a matrix factorization of the product $gt$ using the refined multiplicative tensor product $\overline{\otimes}$. By Lemma \ref{lemma size of X overline tensor Y}, the matrix factors of the product $gt$ will be of size $(4)(4)=16$ since $P$ and $N$ (which are respectively matrix factorizations of $g$ and $t$) are each of size $4$. \\
       We have:\\
$P \overline{\otimes} N = (\phi_{g},\psi_{g})\overline{\otimes}(\phi_{t},\psi_{t})=(\phi_{gt},\psi_{gt})$ where $$(\phi_{gt},\psi_{gt})= (\begin{bmatrix}
    \phi_{g}\otimes \phi_{t}
  \end{bmatrix},\begin{bmatrix}
   \psi_{g}\otimes \psi_{t}
  \end{bmatrix})$$
with \\
$\,\,\,\,\,\,\,\,\,\,\phi_{g}\otimes \phi_{t}\,\,\,\,\,\,\,\,\,\,\,\,\,\,\,\,\,\,=
\begin{bmatrix}
  x & -x^{2} & -y & 0\\
    z &   y^{2}    &  0 & -y\\
 z^{2} &  0    &  y^{2} &  x^{2}\\
     0 &  z^{2}&  -z & x \end{bmatrix} \otimes \begin{bmatrix}
     x^{2} & -y & -y^{2} & 0\\
    y &  z  &  0 & -y^{2}\\
   z &  0  &  z &  y\\
     0 &  z&  -y & x^{2}
  \end{bmatrix}$

  \begin{gather*}
  \setlength{\arraycolsep}{1.0\arraycolsep}
  \renewcommand{\arraystretch}{1.5}
  \text{\footnotesize$\displaystyle
    i.e., \,\phi_{g}\otimes \phi_{t}=\begin{pmatrix}
      x^{3} & -xy & -xy^{2} & 0 & -x^{4} & x^{2}y & x^{2}y^{2} & 0 & -yx^{2} & y^{2} & y^{3} & 0 & 0& 0& 0&0 \\
    xy &   xz    &  0 & -xy^{2} & -x^{2}y & -x^{2}z & 0 & x^{2}y^{2} & -y^{2}& -yz & 0 & y^{3} & 0& 0& 0&0 \\
    xz &  0 &  xz &  xy & -x^{2}z & 0 & -x^{2}z & -x^{2}y & -yz& 0& -yz & -y^{2} & 0& 0& 0&0 \\
     0 &  xz &  -xy & x^{3} & 0 & -x^{2}z & x^{2}y & -x^{4} & 0& -yz& y^{2} & -yx^{2} & 0& 0& 0&0 \\
     zx^{2} & -zy & -zy^{2} & 0 & y^{2}x^{2} & -y^{3} & -y^{4} & 0 & 0& 0& 0 & 0 & -yx^{2}& y^{2}& y^{3}&0 \\
     zy & z^{2} & 0 & -zy^{2} & y^{3} & y^{2}z & 0 & -y^{4} & 0& 0& 0 & 0 & -y^{2}& -yz& 0&y^{3} \\
     z^{2} & 0 & z^{2} & zy & y^{2}z & 0 & y^{2}z & y^{3} & 0& 0& 0 & 0 & -yz& 0& -yz&-y^{2} \\
      0 & z^{2} & -zy & zx^{2} &  0 &  y^{2}z &  -y^{3} & y^{2}x^{2} & 0 & 0& 0& 0&0 & -yz& y^{2} & -yx^{2}\\
      z^{2}x^{2} & -z^{2}y& -z^{2}y^{2} & 0 & 0& 0& 0&0 & y^{2}x^{2} & -y^{3} & -y^{4} & 0 & x^{4} & -x^{2}y & -x^{2}y^{2} & 0\\
     z^{2}y& z^{3}& 0 & -z^{2}y^{2} & 0& 0& 0&0 & y^{3} &   y^{2}z    &  0 & -y^{4} & x^{2}y & x^{2}z & 0 & -x^{2}y^{2} \\
     z^{3}& 0& z^{3} & z^{2}y & 0& 0& 0&0 & y^{2}z &  0 &  y^{2}z &  y^{3} & x^{2}z & 0 & x^{2}z & x^{2}y\\
     0& z^{3}& -z^{2}y & z^{2}x^{2} & 0& 0& 0&0 & 0 &  y^{2}z &  -y^{3} & y^{2}x^{2} & 0 & x^{2}z & -x^{2}y & x^{4}\\
     0& 0& 0 & 0 & z^{2}x^{2} & -z^{2}y& -z^{2}y^{2}&0 & -zx^{2} & zy & zy^{2} & 0 & x^{3} & -xy & -xy^{2} & 0\\
     0& 0& 0 & 0 & z^{2}y& z^{3}& 0&-z^{2}y^{2} &-zy  & -z^{2} & 0 & zy^{2} & xy & xz & 0 & -xy^{2} \\
     0& 0& 0 & 0 & z^{3}& 0& z^{3}&z^{2}y & -z^{2} & 0 & -z^{2} & -zy & xz & 0 & xz & xy \\
      0& 0& 0 & 0 & 0& z^{3}& -z^{2}y& z^{2}x^{2} & 0& -z^{2} & zy & -zx^{2} &  0 &  xz&  -xy & x^{3}
 \end{pmatrix}
  $}
\end{gather*}

  And
  $\,\,\,\,\,\,\,\,\,\,\psi_{g}\otimes \psi_{t}\,\,\,\,\,\,\,\,\,\,\,\,\,\,\,\,\,\,=
  (\begin{bmatrix}
     y^{2} & x^{2} &  y & 0\\
    -z &   x   &  0 & y\\
 -z^{2} &  0    &  x &  -x^{2}\\
     0 &  -z^{2}&  z & y^{2}
  \end{bmatrix}\otimes \begin{bmatrix}
      z & y &  y^{2} & 0\\
    -y &  x^{2}   &  0 & y^{2}\\
 -z &  0    &  x^{2} &  -y\\
     0 &  -z &  y & z
  \end{bmatrix})$

   \begin{gather*}
  \setlength{\arraycolsep}{1.0\arraycolsep}
  \renewcommand{\arraystretch}{1.5}
  \text{\footnotesize$\displaystyle
    i.e., \,\psi_{g}\otimes \psi_{t}=\begin{pmatrix}
       y^{2}z& y^{3} & y^{4} & 0 & x^{2}z & x^{2}y & x^{2}y^{2} & 0 & yz & y^{2} & y^{3} & 0 & 0& 0& 0&0 \\
    -y^{3} &   y^{2}x^{2}  &  0 & y^{4} & -x^{2}y & x^{4} & 0 & x^{2}y^{2} & -y^{2}& yx^{2} & 0 & y^{3} & 0& 0& 0&0 \\
    -y^{2}z &  0 &  y^{2}x^{2} &  -y^{3} & -x^{2}z & 0 & x^{4} & -x^{2}y & -yz& 0& yx^{2} & -y^{2} & 0& 0& 0&0 \\
     0 &  -y^{2}z &  y^{3} & y^{2}z  & 0 & -x^{2}z & x^{2}y & x^{2}z & 0& -yz& y^{2} & yz & 0& 0& 0&0 \\
     -z^{2} & -zy & -zy^{2} & 0 & xz & xy & xy^{2} & 0 & 0& 0& 0 & 0 & yz& y^{2}& y^{3}&0 \\
     zy & -zx^{2} & 0 & -zy^{2} & -xy & x^{3} & 0 & xy^{2} & 0& 0& 0 & 0 & -y^{2}& yx^{2}& 0&y^{3} \\
     z^{2} & 0 & -zx^{2} & zy & -xz & 0 & x^{3} & -xy & 0& 0& 0 & 0 & -yz& 0& yx^{2} &-y^{2} \\
      0 & z^{2} & -zy & -z^{2} &  0 &  -xz &  xy & xz & 0 & 0& 0& 0&0 & -yz& y^{2} & yz\\
      -z^{3} & -z^{2}y& -z^{2}y^{2} & 0 & 0& 0& 0&0 & xz & xy & xy^{2} & 0 & -x^{2}z & -x^{2}y & -x^{2}y^{2} & 0\\
     z^{2}y& -z^{2}x^{2}& 0 & -z^{2}y^{2} & 0& 0& 0&0 & -xy &   x^{3}    &  0 & xy^{2} & x^{2}y & -x^{4} & 0 & -x^{2}y^{2} \\
     z^{3}& 0& -z^{2}x^{2} & z^{2}y & 0& 0& 0&0 & -xz &  0 &  x^{3} &  -xy & x^{2}z & 0 & -x^{4} & x^{2}y\\
     0& z^{3}& -z^{2}y & -z^{3} & 0& 0& 0&0 & 0 &  -xz &  xy & xz & 0 & x^{2}z & -x^{2}y & -x^{2}z\\
     0& 0& 0 & 0 & -z^{3} & -z^{2}y& -z^{2}y^{2}&0 & z^{2} & zy & zy^{2} & 0 & y^{2}z & y^{3} & y^{4} & 0\\
     0& 0& 0 & 0 & z^{2}y& -z^{2}x^{2}& 0&-z^{2}y^{2} &-zy  & zx^{2} & 0 & zy^{2} & -y^{3} & y^{2}x^{2} & 0 & y^{4} \\
     0& 0& 0 & 0 & z^{3}& 0& -z^{2}x^{2} &z^{2}y & -z^{2} & 0 & zx^{2} & -zy & -y^{2}z & 0 & y^{2}x^{2} & -y^{3} \\
      0& 0& 0 & 0 & 0& z^{3}& -z^{2}y& -z^{3} & 0& -z^{2} & zy & z^{2} &  0 &  -y^{2}z &  y^{3} & y^{2}z
 \end{pmatrix}
  $}
\end{gather*}

 \item Next, from the refined algorithm given in the proof of Theorem \ref{thm refined algo for summand-red polyn} , we now need to find a matrix factorization of $r=x^{5}y^{2}$ (which is the first summand in $f$). Evidently, $L=(\phi_{r},\psi_{r})=([x^{5}],[y^{2}])$ is a $1\times 1$ matrix factorization of $x^{5}y^{2}$.
\item Finally, from our algorithm we find a matrix factorization of $f$ by computing $L\widehat{\otimes} (P\overline{\otimes} N)$ which will be of size $2(1)(16)=32$ by Lemma 2.1 of \cite{fomatati2019multiplicative} since $L$ is of size $1$ and $(P\overline{\otimes} N)$ is of size $16$.\\
    We have:
    \begin{align*}
    L\widehat{\otimes} (P\overline{\otimes} N) &= (\phi_{r},\psi_{r})\widehat{\otimes} (\phi_{gt}, \psi_{gs})\\
&=( \begin{bmatrix}
    \phi_{r}\otimes 1_{16}  &  1_{1}\otimes \phi_{gt}      \\
   -1_{1}\otimes \psi_{gt}  &  \psi_{r}\otimes 1_{16}
\end{bmatrix},
\begin{bmatrix}
    \psi_{r}\otimes 1_{16}  &  -1_{1}\otimes \phi_{gt}     \\
    1_{1}\otimes \psi_{gt}  &  \phi_{r}\otimes 1_{16}
\end{bmatrix}
)\\
&=( \begin{bmatrix}
    x^{5}\otimes 1_{16}  &  1\otimes \phi_{gt}      \\
   -1\otimes \psi_{gt}  &  y^{2}\otimes 1_{16}
\end{bmatrix},
\begin{bmatrix}
    y^{2}\otimes 1_{16}  &  -1\otimes \phi_{gt}     \\
    1\otimes \psi_{gt}  &  x^{5}\otimes 1_{16}
\end{bmatrix}
)\\
&=( \begin{bmatrix}
    x^{5}\otimes 1_{16}  &   \phi_{gt}      \\
   - \psi_{gt}  &  y^{2}\otimes 1_{16}
\end{bmatrix},
\begin{bmatrix}
    y^{2}\otimes 1_{16}  &  -\phi_{gt}     \\
    \psi_{gt}  &  x^{5}\otimes 1_{16}
\end{bmatrix}
)\\
&=(\phi_{rgt},\psi_{rgt})
\end{align*}

Where:\\
$\bullet$ $x^{5}\otimes 1_{16}$ (respectively $y^{2}\otimes 1_{16}$) is a $16\times 16$ diagonal matrix with $x^{5}$ (respectively $y^{2}$) on its entire diagonal. \\
$\bullet$ $\phi_{gt}$ and $\psi_{gt}$ were computed above. \\
%
\end{enumerate}
\end{example}
Hence, we found a $32\times 32$ matrix factorization of $f$ viz. a matrix factorization of $f$ of size 32.
%


\begin{quote}
  \textbf{Acknowledgments}
\end{quote}

I am sincerely grateful to the anonymous referee for the insightful comments that help to reshape this paper.



\bibliography{fomatati_ref}
\addcontentsline{toc}{section}
{References}

\end{document}